\documentclass[12pt]{article}

\usepackage{amsmath,amssymb,amsthm}
\usepackage[margin=1.2in]{geometry}
\usepackage{titling}                
\usepackage[osf]{newtxtext} 
\usepackage{zlmtt}          
\usepackage[libertine,bigdelims,cmintegrals,vvarbb]{newtxmath}
\usepackage[cal=euler,frak=euler]{mathalfa} 
\usepackage{bm} 

\usepackage{url}                    
\usepackage{graphicx}               
\usepackage{tikz}                   
\usetikzlibrary{calc,patterns,shadings}

\title{Archimedes' Revenge}

\author{Mark B. Villarino,$^1$ Joseph C. Várilly$^1$
and Runze Li$^2$\\[6pt]
{\normalsize 
$^1\,$Escuela de Matemática, Universidad de Costa Rica,
10101 San José, Costa Rica}\\[3pt]
{\normalsize 
$^2\,$College of Creative Studies, 
University of California at Santa Barbara,}\\
{\normalsize 
Santa Barbara, CA 93106, USA}}

\date{}

\makeatletter
\def\section{\@startsection{section}{1}{\z@}{-3.5ex plus -1ex minus
			  -.2ex}{2.3ex plus .2ex}{\large\bf}}
\def\subsection{\@startsection{subsection}{2}{\z@}{-3.25ex plus -1ex
			  minus -.2ex}{1.5ex plus .2ex}{\normalsize\bf}}
\makeatother

\DeclareMathOperator{\Area}{Area}   
\DeclareMathOperator{\Vol}{Vol}     

\newcommand{\bull}{{\scriptstyle\bullet}} 
\renewcommand{\leq}{\leqslant}      
\newcommand{\otto}{\leftrightarrow} 
\newcommand{\sst}{\scriptstyle}     
\newcommand{\unl}{\underline}       
\newcommand{\x}{\times}             
\renewcommand{\.}{\cdot}            

\bmdefine{\RR}{R}                   

\newcommand{\bR}{\mathbb{R}}        

\newcommand{\half}{{\mathchoice{\thalf}{\thalf}{\shalf}{\shalf}}}
\newcommand{\shalf}{{\scriptstyle\frac{1}{2}}} 
\newcommand{\thalf}{\tfrac{1}{2}}   
\newcommand{\third}{\tfrac{1}{3}}   

\newcommand{\word}[1]{\quad\mbox{#1}\quad} 

\begin{document}

\maketitle

\vspace*{-3pc}

\section{Introduction} 
\label{sec:intro}

Some 2300 years ago the great Archimedes proved that the volume of the
orthogonal intersection of \textit{two} circular cylinders of equal
radius, today called a \textit{bicylinder}, is \textit{two thirds of
the volume of its circumscribed cube}~\cite{Heath60}. (It seems have
remained unnoticed that if the two cylinders intersect \textit{at an
angle}, then the volume of the solid is \textit{two thirds that of the
circumscribed box}. This follows at once by making an affine
transformation of~$\bR^3$, which preserves ratios of volumes.) Later,
modern mathematicians found the volume of the intersection of
\textit{three} such cylinders, a so-called \textit{tricylinder} or
``Steinmetz solid'', and its computation has caused countless
desperation headaches in generations of calculus students.

In 2017, Oliver Knill proposed, under the title ``Archimedes'
Revenge,'' that one prove that \textit{the volume of the intersection
$\RR$ of the three (solid) hyperboloids}
\begin{equation}
x^2 + y^2 - z^2 \leq 1; \quad
y^2 + z^2 - x^2 \leq 1; \quad
z^2 + x^2 - y^2 \leq 1
\label{eq:hyper-trophy} 
\end{equation}
is equal to $\unl{\log 256}$. 

He proposed it as a challenge for the Harvard Maths~21A summer school
in~2017. One of Knill's students (the third author) offered an
integral~\cite{KnillLi17} which computes the volume (see below); this
motivated the solution we give here.

\begin{figure}[htb]
\centering
\includegraphics[width=4.5cm]{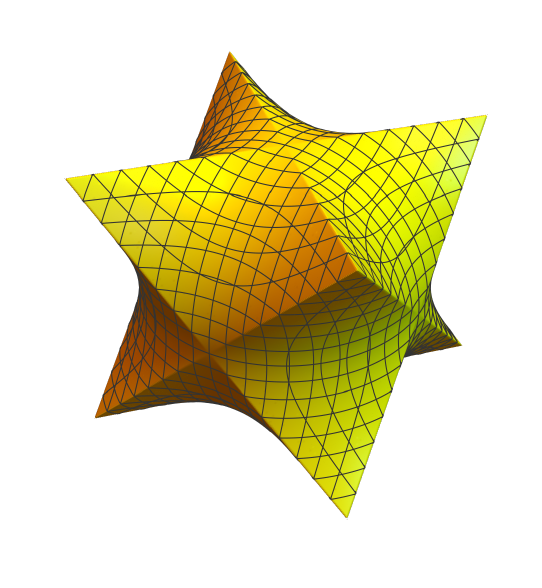}
\caption{The solid contained by the three hyperboloids}
\label{fg:revenge-solid} 
\end{figure}

 
\section{The volume in the first octant} 
\label{sec:first-octant}

By symmetry, the intersection~$\RR$ of the three hyperboloids is the
union of eight congruent solids, one in each octant. We will compute
the volume of the component solid $\RR_1$ in the first octant. The
intricate internal complexity of the solid $\RR_1$ is shown by a
brute-force triple integral for its volume~$V_1$, namely:
\begin{align*}
V_1 &= \int_0^{1/\sqrt{2}} \Biggl\{ \int_0^y \int_0^{\sqrt{x^2+y^2}}
+ \int_y^{1/\sqrt{2}} \int_{\sqrt{y^2-x^2}}^{\sqrt{x^2+y^2}}
+ \int_{1/\sqrt{2}}^{\sqrt{y^2+\half}} 
\int_{\sqrt{x^2-y^2}}^{\sqrt{1-x^2+y^2}} \Biggr\} 1 \,dy\,dx\,dz
\\
&\enspace + \int_{1/\sqrt{2}}^{\sqrt{3/2}} \Biggl\{ 
\int_{\sqrt{y^2-\half}}^{1/\sqrt{2}}
\int_{\sqrt{y^2-x^2}}^{\sqrt{1+x^2-y^2}}
+ \int_{1/\sqrt{2}}^y \int_{\sqrt{x^2+y^2-1}}^{\sqrt{1+x^2-y^2}}
+ \int_y^1 \int_{\sqrt{x^2+y^2-1}}^{\sqrt{1-x^2+y^2}} \Biggr\}
1 \,dy\,dx\,dz
\\
&\enspace + \int_{\sqrt{3/2}}^1 \Biggl\{ 
\int_{\sqrt{y^2-\half}}^{1/\sqrt{2}}
\int_{\sqrt{y^2-x^2}}^{\sqrt{1+x^2-y^2}}
+ \int_{1/\sqrt{2}}^y \int_{\sqrt{x^2+y^2-1}}^{\sqrt{1+x^2-y^2}}
+ \int_y^1 \int_{\sqrt{x^2+y^2-1}}^{\sqrt{1-x^2+y^2}} \Biggr\}
1 \,dy\,dx\,dz.
\end{align*}

This is by no means the end of the story. The first two integrals in
the second and third lines cannot be computed directly, but one must
change the order of integration to calculate them. Thus a complete
computation of the volume $V_1$ by this direct approach is quite
daunting and tedious.


\section{Symmetry} 
\label{sec:symmetry}

We will show how symmetry considerations reduce the computation
of~$V_1$ to a~\textit{single} integral!

The solid $\RR$, whose total volume must be determined, is shown in
Figure~\ref{fg:revenge-solid} -- that appears on the
webpage~\cite{KnillLi17} of Oliver Knill, who kindly contributed the
image. This Mathematica graphic shows its main features: a highly
symmetrical solid whose boundary is a framework with several line
segments supporting negatively curved surface patches taken from the
three hyperboloids. (These line segments also form the edges of
Kepler's \textit{stella octangula}.)

The origin of $\bR^3$ is located at the (hidden) center of the solid.
It evidently has the reflection symmetries $x \otto -x$, 
$y \otto -y$, $z \otto -z$; and the cyclic symmetry 
$(x,y,z) \mapsto (y,z,x)$ of the $120^\circ$ rotation around the 
diagonal line $x = y = z$. 

The key to understanding the solid in Figure~\ref{fg:revenge-solid} is
that a one-sheeted hyperboloid is a \textit{ruled surface} (indeed, a
doubly ruled surface: it can be generated by either of two families of
nonintersecting straight lines). Consider the intersection of any two
of the hyperboloids of Eq.~\eqref{eq:hyper-trophy}:
$$
\biggl\{ \begin{aligned}
y^2 + z^2 - x^2 &= 1 \\
z^2 + x^2 - y^2 &= 1 \end{aligned} \biggr\}
\iff  \biggl\{ \begin{aligned}
z^2 &= 1 \\ x^2 &= y^2 \end{aligned} \biggr\}
\iff  \biggl\{ \begin{gathered}
z = \pm 1 \\ x \pm y = 0 \end{gathered} \biggr\}.
$$
The two lines $z = 1$, $x \pm y = 0$ are the horizontal lines at the
upper boundary in Figure~\ref{fg:revenge-solid}; the other two form
the lower boundary at the bottom. The solid is constrained to lie
inside the third hyperboloid $x^2 + y^2 - z^2 \leq 1$, which cuts off
these four lines at $x^2 + y^2 \leq 2$, yielding four line segments
with eight endpoints $(\pm 1, \pm 1, \pm 1)$, all signs being allowed.

\medskip

It should now be clear that the portion of the solid in the first
octant, $\RR_1$, contains one vertex $(1,1,1)$ of Kepler's star, and
three cross-points of the boundary segments, namely, the standard
basis vectors of~$\bR^3$.

The solid $\RR_1$ is composed of \textit{five}~pieces: 
\begin{itemize}
\item
\textit{two back-to-back tetrahedra} $\Pi_1$ and $\Pi_2$ with a common
base, namely the equilateral triangle with vertices $(1,0,0)$, 
$(0,1,0)$, $(0,0,1)$; and opposite vertices $(0,0,0)$ and $(1,1,1)$,
respectively. Their union $\Pi_1 \cup \Pi_2$ is a triangular 
dipyramid (Figure~\ref{fg:dipyramid}).
\item
\textit{three congruent curved pieces} $S_1$, $S_2$, $S_3$, each of
which is bounded by an outer face of the larger tetrahedron~$\Pi_2$
and one of the hyperboloids.
\end{itemize}

\begin{figure}[htb]
\centering
\begin{tikzpicture}
\coordinate (O) at (0,0) ;
\coordinate (A) at (-0.6,-1) ;
\coordinate (B) at (2.2,0) ;
\coordinate (C) at (0,2) ;
\coordinate (D) at (1.6,1.4) ;
\coordinate (P) at (-0.8,0.8) ;
\coordinate (Q) at (2.5,0.8) ;
\coordinate (X) at (-0.96,-1.6) ;
\coordinate (Y) at (3.1,0) ; 
\coordinate (Z) at (0,2.7) ;
\draw[gray,->] (A) -- (X) node[below left, black] {$x$} ; 
\draw[gray,->] (B) -- (Y) node[right, black] {$y$} ; 
\draw[gray,->] (C) -- (Z) node[above, black] {$z$} ; 
\begin{scope}[every node/.style={font=\footnotesize}]
\draw[thick] (A) node[below right] {$(1,0,0)$}
   -- (B) node[below] {$(0,1,0)$}
   -- (C) node[above right] {$(0,0,1)$}
   -- cycle ; 
\draw[thick,dashed] (O) node[above=3pt, left=-1pt] {$(0,0,0)$} -- (A);
\draw[thick,dashed] (B) -- (O) -- (C) ; 
\draw[thick] (A) -- (D) node[above right] {$(1,1,1)$}
   (B) -- (D) -- (C) ; 
\draw[gray, ->] (P) node[black]{$\Pi_1$} ++(0.16,-0.06)
   -- ++(0.48,-0.18) ;
\draw[gray, ->] (Q) node[black]{$\Pi_2$} ++(-0.24,0) -- ++(-0.6,0) ;
\end{scope}
\foreach \pt in {A,B,C,D,O}  \draw (\pt)  node {$\bull$} ;
\end{tikzpicture} 
\caption{The two adjacent tetrahedra:
         $\Pi_1$ (hidden) and $\Pi_2$ (in the foreground)}
\label{fg:dipyramid} 
\end{figure}

The volumes of the tetrahedra are: $\Vol(\Pi_1) = \dfrac{1}{6}$ and 
$\Vol(\Pi_2) = \dfrac{1}{3}\,$, giving a total of~$\dfrac{1}{2}\,$.


\begin{figure}[tb]
\centering
\begin{tikzpicture}[x=1.4cm, y=0.8cm, z=-0.5cm, scale=2.0]
\coordinate (O) at (0,0,0) ;
\coordinate (A) at (0,0,1) ;
\coordinate (C) at (0,1,0) ;
\coordinate (D) at (1,1,1) ;
\coordinate (F) at (1,0,1) ;
\coordinate (X) at (0,0,1.5) ;
\coordinate (Y) at (0.8,0,0) ; 
\coordinate (Z) at (0,1.4,0) ;
\draw[gray,->] (O) -- (X) ; 
\draw[gray,->] (O) -- (Y) ; 
\draw[gray,->] (O) -- (Z) ; 
\fill[gray!40] (A) -- (F) -- (O) ; 
\draw[gray, dashed] (D) -- (F) ; 
\draw[thick,dashed] (C) -- (A) ;
\draw[thick] (A) -- (D) -- (C) ; 
\draw[thick] plot[variable=\t,smooth,domain=0:90]
    (0,{cos(\t)},{sin(\t)}) ; 
\foreach \pt in {A,C,D}  \draw (\pt)  node {$\bull$} ;
\foreach \ess in {5,10,15,20,25,30,35,40,45,50,55,60,65,70,75,80,85} 
\draw[red] (0,{cos(\ess)},{sin(\ess)}) 
   -- ({sec(\ess) - tan(\ess)},{sec(\ess) - tan(\ess)},1) ;
\foreach \arr in {5,10,15,20,25,30,35,40,45,50,55,60,65,70,75,80,85} 
\draw[blue] (0,{sin(\arr)},{cos(\arr)}) 
   -- ({sec(\arr) - tan(\arr)},1,{sec(\arr) - tan(\arr)}) ;
\end{tikzpicture} 
\caption{The curved piece $S_2$ in the first octant}
\label{fg:curved-roof} 
\end{figure}

\section{Volume of one curved piece} 
\label{sec:curved-piece}

We shall compute the volume of the solid piece $S_2$ bounded by: the
face of~$\Pi_2$ with vertices $(1,0,0)$, $(0,0,1)$, $(1,1,1)$; a
segment of the unit circle in the $xz$-plane; and the hyperboloid
$z^2 + x^2 - y^2 = 1$. The face of $\Pi_2$ lies in the plane 
$x - y + z = 1$. The projection of this tetrahedral face onto the
$xy$-plane is the triangle with vertices $(0,0)$, $(1,0)$, $(1,1)$.
See Figure~\ref{fg:curved-roof}.

The portion of the hyperboloid that forms the roof of~$S_2$ is shown 
in Figure~\ref{fg:curved-roof} with both rulings by line segments. 
Its curved boundary (a quarter circle) can be parametrized by either
$(\sin\theta,0,\cos\theta)$ or $(\cos\phi,0,\sin\phi)$ with the 
angles in the interval $[0,\pi/2]$. The rulings proceed from there to 
either of the straight sides, as follows:
\begin{alignat*}{2}
t &\mapsto (\sin\theta + t\cos\theta, 0, \cos\theta - t\sin\theta),
&& 0 \leq t \leq \sec\theta - \tan\theta;
\\
s &\mapsto (\cos\phi - s\sin\phi, 0, \sin\phi + s\cos\phi),
& \qquad & 0 \leq s \leq \sec\phi - \tan\phi.
\end{alignat*}

\medskip

The triple integral for the volume of the curved piece $S_2$ is thus:
$$
\boxed{I := \int_0^1 \int_y^1 \int_{1-x+y}^{\sqrt{1+y^2-x^2}}
1 \,dz\,dx\,dy.}
$$
The inner integration is immediate:
$$
I = \int_0^1 \int_y^1 \sqrt{1 + y^2 - x^2} - (1 - x + y) \,dx\,dy.
$$
We write
$$
I_1 := \int_0^1 \int_y^1 \sqrt{1 + y^2 - x^2} \,dx\,dy
\word{and}
I_2 := \int_0^1 \int_y^1 (1 - x + y) \,dx\,dy.
$$
A routine calculation gives $I_2 = \dfrac{1}{3}\,$.

Using the indefinite integral of $\sqrt{a^2 - x^2}$ with
$a = \sqrt{1 + y^2}$ and then integrating by parts twice, we obtain
\begin{align*}
I_1 &= \int_0^1 \biggl[ \frac{x}{2} \sqrt{1 + y^2 - x^2}
+ \frac{1 + y^2}{2} \arcsin\biggl( \frac{x}{\sqrt{1 + y^2}} \biggr)
\biggr]_{x=y}^{x=1} \,dy
\\
&= \int_0^1 \biggl\{ \frac{1 + y^2}{2}
\arcsin\biggl( \frac{1}{\sqrt{1 + y^2}} \biggr)
- \frac{1 + y^2}{2} \arcsin\biggl( \frac{y}{\sqrt{1 + y^2}} \biggr)
\biggr\} \,dy
\\
&= \frac{\pi}{6} + \frac{1}{2} \int_0^1
\biggl(\frac{y + \third y^3}{1 + y^2} \biggr) \,dy
- \biggl\{ \frac{\pi}{6} - \frac{1}{2} \int_0^1 
\biggl(\frac{y + \third y^3}{1 + y^2} \biggr) \,dy \biggr\}
\\
&= \int_0^1 \biggl(\frac{y + \third y^3}{1 + y^2} \biggr) \,dy
= \frac{1}{3} \int_0^1 \biggl( \frac{2y}{1 + y^2} + y \biggr) \,dy 
= \frac{\log 2}{3} + \frac{1}{6}\,.
\end{align*}

Thus the volume of the curved piece $S_2$ is
$$
I = I_1 - I_2 = \frac{\log 2}{3} + \frac{1}{6} - \frac{1}{3}
=  \frac{\log 2}{3} - \frac{1}{6}\,.
$$


\section{Total Volume} 
\label{sec:total-volume}

Therefore the volume of the solid~$\RR_1$ in the first quadrant is
$$
V_1 = 3I + \frac{1}{2} 
= 3\biggl( \frac{\log 2}{3} - \frac{1}{6} \biggr) + \frac{1}{2} 
= \log 2;
$$
and the total volume of the solid~$\RR$ in Archimedes' revenge is:
$$
\Vol(\RR) = 8\,V_1 = 8 \log 2 = \log 256.
\eqno \text{\textbf{qed!}}
$$


\section{Comment} 
\label{sec:final-comment}

The idea to exploit the symmetry of a component in one octant was
already suggested by the third author~\cite{KnillLi17}. His solution
stated that the solid~$\RR_1$ is composed of a tetrahedron (actually a
dipyramid) of volume~$\half$ and the three congruent curved pieces. He
computed the volume of a curved piece using the following integral:
$$
I := \frac{1}{2} \int_0^1 \biggl\{ (z^2 + 1) \biggl(
\frac{\pi}{2} - 2\arctan z \biggr) + z^2 - 1 \biggr\} \,dz
\eqno (*)
$$
which the first two authors found somewhat mysterious, and indeed
their attempt to decipher it led to the solution offered here.


\section*{Acknowledgements}

We are grateful to Oliver Knill for sharing the image in
Figure~\ref{fg:revenge-solid} that appears on his
webpage~\cite{KnillLi17} (which also includes the student's suggested
solution), and for helpful comments. We also acknowledge support from
the Vicerrectoría de Investigación of the Universidad de Costa~Rica.



\section*{Postscript}

A shortened version of the above solution was published as a Classroom
Note by the first two authors in the \textit{College Mathematical
Journal}, vol.~56 (2024), 257--259. Shortly afterward, in an email
exchange, the third author sent them his original solution, employing
the aforementioned integral~$(*)$. We append this elegant solution
below.

\newpage

\appendix

\section{The original calculation} 
\label{app:early-solution}

Firstly, we can divide the solid into $8$ parts, one part is shown 
in Figure~\ref{fg:first-octant}:

\begin{figure}[htb]
\centering
\begin{tikzpicture}[x=1.3cm, y=1.1cm, z=-0.8cm, scale=2.0]
\coordinate (O) at (0,0,0) ;
\coordinate (A) at (0,0,1) ;
\coordinate (B) at (1,0,0) ;
\coordinate (C) at (0,1,0) ;
\coordinate (D) at (1,1,1) ;
\coordinate (E) at (0,1,1) ;
\coordinate (F) at (1,0,1) ;
\coordinate (G) at (1,1,0) ;
\coordinate (X) at (0,0,1.3) ;
\coordinate (Y) at (1.3,0,0) ; 
\coordinate (Z) at (0,1.3,0) ;
\draw[gray,->] (A) -- (X) node[below left, black] {$x$} ; 
\draw[gray,->] (B) -- (Y) node[right, black] {$y$} ; 
\draw[gray,->] (C) -- (Z) node[above, black] {$z$} ; 
\begin{scope}[every node/.style={font=\footnotesize}]
\draw[thick] (A) node[below right] {$(1,0,0)$}
   -- (B) node[below=2pt] {$(0,1,0)$}
   -- (C) node[above right] {$(0,0,1)$}
   -- cycle ; 
\draw[thick,dashed] (O) node[above=3pt, left=-1pt] {$(0,0,0)$} -- (A);
\draw[thick,dashed] (B) -- (O) -- (C) ; 
\fill[gray!30] (A) -- (B) -- (D) -- cycle ;
\draw[thick] (A) -- (D) node[above right=-2pt] {$(1,1,1)$}
   (B) -- (D) -- (C) ; 
\end{scope}
\draw[very thin, blue] (C) -- (E) -- (A) -- (F) -- (B) -- (G)
   -- cycle  (D) -- (E)  (D) -- (F)  (D) -- (G) ; 
\draw[thick] plot[variable=\t,smooth,domain=0:90]
    (0,{cos(\t)},{sin(\t)}) ; 
\draw[thick] plot[variable=\t,smooth,domain=0:90]
    ({sin(\t)},0,{cos(\t)}) ; 
\draw[thick] plot[variable=\t,smooth,domain=0:90]
    ({cos(\t)},{sin(\t)},0) ; 
\foreach \ess in {5,10,15,20,25,30,35,40,45,50,55,60,65,70,75,80,85} 
\draw[red] ({sin(\ess)},0,{cos(\ess)}) 
   -- (1,{sec(\ess) - tan(\ess)},{sec(\ess) - tan(\ess)}) ;
\foreach \arr in {5,10,15,20,25,30,35,40,45,50,55,60,65,70,75,80,85} 
\draw[blue] ({cos(\arr)},0,{sin(\arr)}) 
   -- ({sec(\arr) - tan(\arr)},{sec(\arr) - tan(\arr)},1) ;
\foreach \pt in {A,B,C,D,O}  \draw (\pt)  node {$\bull$} ;
\end{tikzpicture} 
\caption{The curved piece $S_1$ in the first octant}
\label{fg:first-octant} 
\end{figure}

Next, we need to calculate the volume of the dark part in the picture.
(This is $S_1$, one of the curved pieces mentioned above.)

\medskip

If we slice it with a horizontal plane (at height~$z$ between $0$
and~$1$), we can get a shape like the one below
(Figure~\ref{fg:sliver}). Call $S(z)$ the area of the shaded shape.

\begin{figure}[htb]
\centering
\begin{tikzpicture}[scale=3.0]
\coordinate (A) at (0,0) ;
\coordinate (Bb) at (-66:1.2cm) ;
\coordinate (Cc) at (-24:1.2cm) ;
\coordinate (D) at (1,-1) ;
\coordinate (X) at (0,-1) ;
\coordinate (Y) at (1,0) ; 
\coordinate (B) at (intersection of A--Bb and X--D) ;
\coordinate (C) at (intersection of A--Cc and Y--D) ;
\coordinate (Ab) at ($ (A)!0.6!(B) $) ; 
\coordinate (Ax) at ($ (A)!0.5!(X) $) ;
\coordinate (Ay) at ($ (A)!0.5!(Y) $) ;
\coordinate (Bx) at ($ (B)!0.5!(X) $) ;
\coordinate (T1) at (-12:0.3cm) ;
\coordinate (T2) at (-45:0.25cm) ;
\coordinate (T3) at (-78:0.3cm) ;
\draw (A) node[above left] {$A$} -- (X) -- (D) -- (Y) -- cycle ;
\draw (B) node[below] {$B$} -- (A) -- (C) node[right] {$C$} ;
\draw[thick, pattern={vertical lines}, pattern color=blue]
   (B) arc(-66:-24:1.095cm) -- cycle ; 
\draw[very thin, blue] (A) ++(-90:0.2cm) arc(-90:-66:0.2cm) ;
\draw[thin] (A) ++(-66:0.15cm) arc(-66:-24:0.15cm) ;
\draw[very thin, blue] (A) ++(-24:0.2cm) arc(-24:0:0.2cm) ;
\draw (Ab) node[above right=-3pt] {$\sst\sqrt{1 + z^2}$} ; 
\draw (Ax) node[left] {$\sst 1$} ; 
\draw (Ay) node[above] {$\sst 1$} ; 
\draw (Bx) node[above] {$\sst z$} ; 
\draw (T1) node {$\sst\theta$} ; \draw (T3) node {$\sst\theta$} ; 
\draw (T2) node {$\sst\alpha$} ; 
\foreach \pt in {A,B,C,D,X,Y}  \draw (\pt)  node {$\bull$} ;
\end{tikzpicture} 
\caption{A horizontal slice of the curved piece $S_2$}
\label{fg:sliver} 
\end{figure}

Here $A = (0,0,z)$, $B = (1,z,z)$, $C = (z,1,z)$. The parameters in 
Figure~\ref{fg:sliver} are related by
$$
z = \tan\theta, \quad  
\alpha = \tfrac{\pi}{2} - 2\theta = \tfrac{\pi}{2} - 2\arctan z,
$$
and for convenience we put $r := \sqrt{1 + z^2}$, the radius of the 
circular arc $BC$ (which is a slice of the hyperboloid 
$x^2 + y^2 - z^2 = 1$). The area $S(z)$ is that of the circle sector
$ABC$ minus the area of the triangle $\triangle\,ABC$. Since
$\overrightarrow{AB} = (1,z,0)$ and $\overrightarrow{AC} = (z,1,0)$,
the triangle has area 
$$
\Area(\triangle\,ABC) 
= \half\,\bigl\| \overrightarrow{AB} \x \overrightarrow{AC} \bigr\|
= \half(1 - z^2).
$$
The area of the circle sector $ABC$ is $\half r^2 \alpha$, and so
$$
S(z) = \frac{1}{2} \Bigl(
(z^2 + 1)\bigl( \tfrac{\pi}{2} - 2\arctan z \bigr) + (z^2 - 1) \Bigr).
$$

The volume of the curved piece $S_2$ is obtained by integrating $S(z)$
from $z = 0$ to $z = 1$:
\begin{align*}
\Vol(S_2) &= \int_0^1 S(z) \,dz
= \frac{1}{2} \int_0^1 \Bigl( \frac{\pi}{2}(z^2 + 1)
+ (z^2 - 1) - 2(z^2 + 1) \arctan z \Bigr) \,dz
\\
&= \frac{1}{2} \. \frac{\pi}{2} \. 
\biggl[ \frac{z^3}{3} + z \biggr]_{z=0}^{z=1}
+ \frac{1}{2}\, \biggl[ \frac{z^3}{3} - z \biggr]_{z=0}^{z=1}
- \int_0^1 \bigl( (z^2 + 1) \arctan z \bigr) \,dz
\\
&= \frac{\pi-1}{3} - \biggl[
\frac{z^3}{3} \arctan z + z \arctan z \biggr]_{z=0}^{z=1}
+ \int_0^1 \Bigl( \frac{z^3}{3} + z \Bigr) \,d(\arctan z)
\\
&= \frac{\pi-1}{3} - \frac{\pi}{12} - \frac{\pi}{4}
+ \int_0^1 \frac{\frac{1}{3} z^3 + z}{z^2 + 1} \,dz
= -\frac{1}{3} + \frac{1}{3} \int_0^1 
\Bigl( z + \frac{2z}{z^2 + 1} \Bigr) \,dz
\\
&= -\frac{1}{3} + \frac{1}{3} \biggl[
\frac{z^2}{2} + \log(z^2 + 1) \biggr]_{z=0}^{z=1}
= -\frac{1}{6} + \frac{\log 2}{3}\,.
\end{align*}

\begin{figure}[htb]
\centering
\begin{tikzpicture}[x=1.3cm, y=1.1cm, z=-0.8cm, scale=2.0]
\coordinate (O) at (0,0,0) ;
\coordinate (A) at (0,0,1) ;
\coordinate (B) at (1,0,0) ;
\coordinate (C) at (0,1,0) ;
\coordinate (D) at (1,1,1) ;
\coordinate (E) at (0,1,1) ;
\coordinate (F) at (1,0,1) ;
\coordinate (G) at (1,1,0) ;
\coordinate (X) at (0,0,1.3) ;
\coordinate (Y) at (1.3,0,0) ; 
\coordinate (Z) at (0,1.3,0) ;
\pattern[pattern={horizontal lines}, pattern color=blue!40]
   (C) -- (A) -- (D) -- cycle ;
\pattern[pattern={vertical lines}, pattern color=blue!40]
   (C) -- (E) -- (D) -- cycle ;
\pattern[pattern={north east lines}, pattern color=blue!40]
   (E) -- (A) -- (D) -- cycle ;
\draw[gray,->] (A) -- (X) node[below left, black] {$x$} ; 
\draw[gray,->] (B) -- (Y) node[right, black] {$y$} ; 
\draw[gray,->] (C) -- (Z) node[above, black] {$z$} ; 
\begin{scope}[every node/.style={font=\footnotesize}]
\draw[thick] (A) node[below right] {$(1,0,0)$}
   -- (B) node[below=2pt] {$(0,1,0)$}
   -- (C) node[above right] {$(0,0,1)$}
   -- cycle ; 
\draw[thick,dashed] (O) -- (A);
\draw[thick,dashed] (B) -- (O) -- (C) ; 
\draw[thick] (A) -- (D) node[above right=-2pt] {$(1,1,1)$}
   (B) -- (D) -- (C) ; 
\draw (E) node[left=-1pt] {$(1,0,1)$} ;
\end{scope}
\draw[very thin, blue] (C) -- (E) -- (A) -- (F) -- (B) -- (G)
   -- cycle  (D) -- (E)  (D) -- (F)  (D) -- (G) ; 
\foreach \pt in {A,B,C,D,E,F,G,O}  \draw (\pt)  node {$\bull$} ;
\end{tikzpicture} 
\caption{The unit cube with three equal pyramids removed}
\label{fg:new-dipyramid} 
\end{figure}

Then we can calculate the volume $V$ of the remaining part in the
first octant, which is the cube of side~$1$ with three equal pyramids
removed (Figure~\ref{fg:new-dipyramid}):
$$
V = 1 \x 1 \x 1 - 3 \bigl( \tfrac{1}{3}(\half \x 1 \x 1) \bigr)
= \frac{1}{2} \,.
$$

So the total volume of the part in the first octant (see
Figure~\ref{fg:first-octant}) is:
$$
V + 3\biggl( -\frac{1}{6} + \frac{\log 2}{3} \biggr)
= \frac{1}{2} - \frac{1}{2} + \log 2 = \log 2.
$$
Finally, the total volume of the whole object is
$$
8 \log 2 = \log(2^8) = \log 256.
\eqno \text{\textbf{qed!}}
$$


\begin{thebibliography}{26}

\bibitem{Heath60}
Archimedes,
\textit{The Works of Archimedes: The Method},
T.~L. Heath, ed.,
Dover Publications, New York, 1960.

\bibitem{KnillLi17}
From the webpages of Oliver Knill, 2017:
\url{http://people.math.harvard.edu/~knill/teaching/summer2017/exhibits/revenge/}


\end{thebibliography}
\end{document}